\documentclass[11pt]{article}
\usepackage{amssymb}
\def\bbR{\mathbb{R}}
\def\LTR{L^2(\bbR)}

\def\iint{\mathop{\int\!\!\!\int}}

\newtheorem{Theorem}{Theorem}[section]

\newtheorem{Proposition}[Theorem]{Proposition}
\newtheorem{Lemma}[Theorem]{Lemma}

\newtheorem{Remark}[Theorem]{Remark}

\def\supp{\mathrm{supp}\,}
\makeatletter

       \@addtoreset{equation}{section}
 \makeatother
\begin{document}

\title{Inversion Formula for the  Windowed Fourier Transform\thanks{This work was supported partially by the
National Natural Science Foundation of China(10971105 and 10990012) and the
Natural Science Foundation   of Tianjin (09JCYBJC01000).}}

\author{Wenchang Sun\\
  Department of Mathematics and LPMC,  Nankai University,
      Tianjin~300071, China\\
Email:\,      sunwch@nankai.edu.cn}

\date{ }

\maketitle

\begin{abstract}
In this paper, we study the inversion formula for recovering a function from its windowed Fourier transform.
We give a rigorous proof for an inversion formula which is known in engineering.
We show that the integral involved in the formula is convergent almost everywhere
on $\bbR$  as well as  in $L^p$ for all $1<p<\infty$ if the function to be reconstructed is.
\end{abstract}

\textbf{Keywords}.\,\,
Fourier transforms; windowed Fourier transforms; inversion formula.

\textbf{2000 Mathematics Subject Classification}. \,\,  42A38.

\section{Introduction and the Main Result}
The  Fourier transform  is a very useful mathematical tool, which has been widely used in
characterization of function spaces as well as in signal and image processing
\cite{Duo,Stein}. For a function $f\in L^1(\bbR)$,   the Fourier transform of $f$
is defined by
\[
   \hat{f}(\omega) = \int_{\bbR} f(x) e^{-ix\omega} dx.
\]
To study local properties of functions (signals), the windowed Fourier transform, also known as short-time Fourier
transform, is introduced.

Given a window function $g(x)$,  the windowed Fourier transform of a function $f$ with respect to $g$ is defined by
\[
  (F_gf)(t,\omega) = \int_{\bbR} f(x) \overline{g(x-t)} e^{-ix\omega}dx.
\]
It is easy to see that $F_gf$ is well defined if $f\in L^p(\bbR)$ and $g\in L^{p'}(\bbR)$, where $p$, $p' \ge 1$ and $1/p + 1/p'=1$.

Continuous and discrete   windowed Fourier transforms  have been discussed extensively in the literature since
they are widely used
in communication theory, quantum mechanics, and many other fields.
We refer to \cite{C,Chui,D1,FS,Gro} for an introduction to the windowed Fourier transform.

Finding a computationally efficient algorithm for the inversion of  windowed Fourier transforms
is a fundamental topic in both theory and applications.
The classical method to recover $f$ from its windowed Fourier transform is
to use the following inversion formula,
\begin{equation}\label{eq:r1}
    f(x) = \frac{1}{2\pi\|g\|_2^2} \iint_{\bbR^2} (F_gf)(t,\omega) g(x-t) e^{ix\omega} dtd\omega,
\end{equation}
where we assume that $g\in \LTR$. It can be shown that the convergence is in $L^2(\bbR)$ as well
as in many other spaces if the function to be reconstructed is and  $g$ satisfies some further conditions~\cite{FS}.

Since a double integral is involved in (\ref{eq:r1}), it is obviously very complicated.
An alternate method is to use the filter-bank summation \cite{AR},
\begin{equation}\label{eq:fbs}
    f(x) = \frac{1}{2\pi\overline{g(0)}} \int_{\bbR} (F_gf)(x, \omega) e^{ix\omega} d\omega,
\end{equation}
where we assume that $g(0)\ne 0$.
Note that (\ref{eq:fbs}) was presented in \cite{AR} in a discrete version for compactly supported window functions
and the authors stated that their results may be equally well stated in a continuous
time-domain setting.

Although (\ref{eq:fbs}) is well known in engineering, the convergence of the integral is not well
stated in literature.
In this paper, we show that the integral in (\ref{eq:fbs}) is convergent
 in $L^p(\bbR)$ for all $1<p<\infty$ if the function $f$ is.
Moreover, by applying the Carleson-Hunt theorem, we also show that
the convergence is almost everywhere on $\bbR$.

Before stating our result, we introduce some definitions.
Throughout this paper, $x_0$ is a fixed real number.
For any  $A_1, A_2>0$, define
\begin{equation}\label{eq:t1}
  (T_{A_1,A_2}f)(x) = \int_{-A_1}^{A_2} (F_gf)(x-x_0, \omega) e^{ix\omega} d\omega.
\end{equation}

Our main result is the following.

\begin{Theorem}\label{thm:main}
Suppose that $g$ is continuous and that  $g, \hat g\in L^1(\bbR)$.
Then for any $f\in L^p(\bbR)$, $1<p<\infty$, we have
 \begin{equation}\label{eq:ee1}
 \lim_{ A_1,A_2\rightarrow \infty}
    \|T_{A_1,A_2}f - 2\pi\overline{g(x_0)} f\|_p = 0
 \end{equation}
and
 \begin{equation}\label{eq:ee1a}
 \lim_{ A\rightarrow \infty}
     (T_A f)(x)  =  2\pi\overline{g(x_0)} f(x),\qquad a.e.,
 \end{equation}
where we use the shortcut $T_A f = T_{A,A}f$.
\end{Theorem}

\begin{Remark}
\label{rem:r2}
The reconstruction  formula (\ref{eq:t1}) is stable in the sense that for any $f, \tilde f\in L^p(\bbR)$,
\[
  \| T_{A_1, A_2} (f-\tilde f)\|_{p} \le 2 C_p \|\hat g\|_1\|f-\tilde f\|_{p}, \qquad \forall A_1, A_2>0,
\]
where $C_p$ is a constant depending only on $p$.
 For details, see the proof of Theorem~\ref{thm:main}.
\end{Remark}

In Section 2, we give the proof of Theorem\ref{thm:main}, which is based on the famous Carleson-Hunt theorem \cite{Carl,Hunt}
for Fourier series and the extension to Fourier integrals by Kenig and Tomas \cite{KT}.
\section{Proof of the Main Result}
In this section, we give the proof of the main result.

We begin with a simple lemma on the Fourier transform, for which we omit the proof.
\begin{Lemma} \label{Lm:L0}
For any $f\in\LTR$ with $\hat f\in L^1(\bbR)$, we have
\[
   f(x)=\frac{1}{2\pi}\int_{\bbR} \hat f(\omega) e^{ix\omega} d\omega, \qquad a.e.
\]
\end{Lemma}

We also need the following formula on the windowed Fourier transform.

\begin{Proposition}[{\cite[Lemma 3.1.1]{Gro}}]\label{prop:p2}
For any $f,g\in\LTR$, we have
\[
   (F_gf)(x,\omega) = \frac{1}{2\pi}(F_{\hat g} \hat f)(\omega, -x) e^{-ix\omega}.
\]
\end{Proposition}

Next,  we show that for $f\in \LTR$ with $\hat{f}\in L^1(\bbR)$,
$T_{A_1,A_2}f$ is convergent in $L^{\infty}$ norm.

\begin{Lemma}
\label{Lm:L3}
Suppose that $g$ is continuous and that $g, \hat g\in L^1(\bbR)$.
Then for any  $f\in\LTR$ with $\hat{f}\in L^1(\bbR)$,  we have
 \begin{equation}\label{eq:ee2}
   \lim_{ A_1,A_2 \rightarrow \infty}
    \|T_{A_1,A_2}f - 2\pi\overline{g(x_0)} f\|_{\infty} = 0.
 \end{equation}
\end{Lemma}

\begin{proof}
For any $f\in\LTR$, we see from Proposition~\ref{prop:p2} that
\begin{eqnarray}
(T_{A_1,A_2}f)(x)&=& \int_{-A_1}^{A_2} (F_gf)(x-x_0, \omega) e^{ix\omega} d\omega
   \nonumber \\
&=& \frac{1}{2\pi}\int_{-A_1}^{A_2} (F_{\hat{g}} \hat{f})(\omega, x_0-x) e^{-i(x-x_0)\omega}
   e^{ix\omega}  d\omega
   \nonumber \\
&=& \frac{1}{2\pi}\int_{-A_1}^{A_2} e^{ix_0\omega}  d\omega
   \int_{\bbR}
     \hat{f}(y) \overline{\hat{g}(y-\omega)}  e^{-iy(x_0-x)} dy
    \nonumber \\
&=& \frac{1}{2\pi}
   \int_{\bbR}
     \hat{f}(y)   e^{iyx} dy
     \int_{-A_1}^{A_2}\overline{\hat{g}(y-\omega)} e^{-ix_0(y- \omega)}  d\omega
       \nonumber \\
 &=& \frac{1}{2\pi}
   \int_{\bbR}
     \hat{f}(y)   e^{iyx} dy
     \int_{y-A_2}^{y+A_1} \overline{\hat{g}(\omega)} e^{-ix_0 \omega}  d\omega
    \nonumber \\
&=& \frac{1}{2\pi}
   \int_{\bbR}
 \overline{\hat{g}(\omega)} e^{-ix_0 \omega}  d\omega
     \int_{\omega-A_1}^{\omega+A_2}
          \hat{f}(y)   e^{iyx} dy,
  \label{eq:t2a}
\end{eqnarray}
where we use Fubini's theorem twice.
By Lemma~\ref{Lm:L0},   for almost every $x$,
\begin{eqnarray}
 && (T_{A_1,A_2}f)(x)  -2\pi \overline{g(x_0)} f(x) \nonumber \\
  &=&
   \int_{\bbR}
 \overline{\hat{g}(\omega)} e^{-ix_0 \omega}  d\omega
    \bigg(\frac{1}{2\pi} \int_{\omega-A_1}^{\omega+A_2}
          \hat{f}(y)   e^{iyx} dy -f(x) \bigg) \nonumber \\
  &=&
       \frac{1}{2\pi} \int_{\bbR}
 \overline{\hat{g}(\omega)} e^{-ix_0 \omega}  d\omega
 \int_{{y<\omega-A_1}\atop{\mathrm{or\,}y>\omega+A_2}}
          \hat{f}(y)   e^{iyx} dy.\nonumber
\end{eqnarray}
Hence
\[
  \left\|T_{A_1,A_2}f  -2\pi  g(x_0) f\right\|_{\infty}
  \le
  \frac{1}{2\pi}  \int_{\bbR}
  |\hat{g}(\omega)|  d\omega
 \int_{{y<\omega-A_1}\atop{\mathrm{or\,}y>\omega+A_2}}
          |\hat{f}(y)  | dy.
\]
By the dominated convergence theorem, we get
\[
   \lim_{ A_1, A_2 \rightarrow \infty}
    \|T_{A_1,A_2}f - 2\pi\overline{g(x_0)} f\|_{\infty} = 0.
\]
This completes the proof.
\end{proof}

In the followings we prove the convergence in $L^p(\bbR)$.
First, we  show  that $T_{A_1,A_2}$ is well defined on $L^p(\bbR)$.

\begin{Lemma} \label{Lm:L1}
Suppose that $g$ is continuous and that $g, \hat g\in L^1(\bbR)$.
For any $A_1, A_2>0$, let
\begin{eqnarray}
  K_{A_1,A_2}(x,y) &=& \overline{g(y-x+x_0)}\cdot\bigg(
  \frac{\sin\, A_1(y-x) + \sin\, A_2(y-x)}{(y-x)}\nonumber \\
&&   - i \cdot\frac{2\sin\frac{(A_2-A_1)(y-x)}{2} \sin\frac{(A_2+A_1)(y-x)}{2}}{(y-x)}
    \bigg). \label{eq:KA}
\end{eqnarray}
Then we have
\[
   ( T_{A_1,A_2} f )(x) = \int_{\bbR}
   f(y) K_{A_1,A_2}(x,y) dy,\qquad \forall f\in L^p(\bbR).
\]
\end{Lemma}

\begin{proof}
Since $g, \hat g\in L^1(\bbR)$, we have $g\in L^p(\bbR)$ for all $1<p<\infty$.
Hence $F_gf$ is well defined for any $f\in L^p(\bbR)$. We have
\begin{eqnarray}
   ( T_{A_1,A_2} f )(x)
&=& \int_{-A_1}^{A_2} (F_gf)(x-x_0, \omega) e^{ix\omega} d\omega
        \nonumber \\
&=& \int_{-A_1}^{A_2}   d\omega
   \int_{\bbR} f(y) \overline{g(y-x+x_0)} e^{-iy\omega}  e^{ix\omega} dy
        \nonumber \\
&=& \int_{\bbR}   dy \int_{-A_1}^{A_2}
   f(y) \overline{g(y-x+x_0)} e^{-i(y-x)\omega}    d\omega
        \nonumber \\
&=& \int_{\bbR}
   f(y)  K_{A_1,A_2}(x,y)  dy,
        \nonumber
\end{eqnarray}
where Fubini's theorem is used. This completes the proof.
\end{proof}

The pointwise convergence of Fourier series  is a deep result in harmonic analysis.
Carleson proved that the Fourier series of a function in $L^2[-\pi, \pi]$
is convergent almost everywhere \cite{Carl}. Hunt \cite{Hunt}
extended this result   to  $L^p[-\pi, \pi]$  for $1<p<\infty$.
And Kenig and Tomas \cite{KT} proved the pointwise convergence of Fourier integral on $L^p(\bbR)$.
For our purpose, we cite the Carleson-Hunt theorem in the following form.
\begin{Proposition}\label{prop:CH}
For $A>0$ and $1<p<\infty$, define
\begin{equation}\label{eq:sa}
 ( S_A f)(x) = \int_{\bbR} f(y) \frac{\sin A(x-y)}{\pi(x-y)}dy,\qquad f\in L^p(\bbR).
\end{equation}
Then $S_A$ is a bounded linear operator on $L^p(\bbR)$ and there exists some constant $C_p$ such that
\[
    \Big\| \sup_{A>0} |( S_A f)(x)|  \Big\|_p \le C_p \|f\|_p.
\]
\end{Proposition}

The Fourier multiplier is a useful tool in the study of Fourier transform.
The following result on the Fourier multiplier is useful in studying the
convergence of $T_{A_1, A_2}$.

 \begin{Proposition}[{\cite[Corollary 3.8]{Duo}}] \label{prop:p4}
Suppose that  $h$ is a function of bounded variation on $\bbR$
and that $(Tf)\hat{} = h\cdot \hat f$ for $f\in\LTR$. Then $T$
can be extended to an operator on $L^p(\bbR)$, $1<p<\infty$ and
\[
   \|Tf\|_p \le C_p V_h \|f\|_p,\qquad \forall f\in L^p(\bbR),
\]
where $V_h$ is the total  variation of $h$ on $\bbR$ and
$C_p$ is a constant depending only on $p$.
\end{Proposition}

The following lemma shows that $T_{A_1,A_2} f$ converges to $f$ in $L^p(\bbR)$
whenever $f$ is in $C_c^1(\bbR)$, the space of all continuous
differentiable functions which are compactly
supported.
\begin{Lemma}\label{Lm:L2}
For any $f\in C_c^1(\bbR)$, we have
\begin{equation}\label{eq:e10}
    \lim_{A_1,A_2\rightarrow \infty } \| T_{A_1,A_2} f - 2\pi\overline{g(x_0)} f  \|_p=0,
    \qquad 1<p\le \infty.
\end{equation}
\end{Lemma}

\begin{proof}
Fix some $f\in C_c^1(\bbR)$. Suppose that $\supp f\subset [-\Omega, \Omega]$ for some constant $\Omega>0$.
Since $f, f'\in\LTR$, we have $\hat{f}\in L^1(\bbR)$.
By Lemma~\ref{Lm:L3},
\begin{equation}\label{eq:con1}
  \lim_{A_1,A_2\rightarrow \infty }  \|T_{A_1,A_2}  f - 2\pi\overline{g(x_0)} f \|_{\infty}=0.
\end{equation}

Next we assume that $1<p<\infty$. By (\ref{eq:con1}), we have
\begin{equation}\label{eq:s1}
   \lim_{A_1,A_2\rightarrow \infty } \|(T_{A_1,A_2} f - 2\pi\overline{g(x_0)} f )\cdot   \chi^{}_{[-2\Omega, 2\Omega]}\|_p=0.
\end{equation}
On the other hand,  put
\[
  K(x,y) = \frac{4\|g\|_{\infty}}{|x-y|},\qquad x\ne y.
\]
By   Minkovski's inequality, we have
\begin{eqnarray*}
&&\left(\int_{|x|\ge 2\Omega} \left| \int_{|y|\le \Omega} K (x,y) |f(y)|dy\right|^p dx\right)^{1/p}\\
&\le&
  \int_{|y|\le \Omega} |f(y)| \left(\int_{|x|\ge 2\Omega} | K (x,y)|^p dx \right)^{1/p}dy\\
&\le& 4\|g\|_{\infty}  \int_{|y|\le \Omega} |f(y)|  \left( \int_{|x|>2\Omega} \frac{dx}{|x-y|^p}\right)^{1/p}dy \\
&=& M_p \|f\|_1\\
&\le& M_p (2\Omega)^{1/p'}\|f\|_p,
\end{eqnarray*}
where $M_p$ is a constant and $1/p+1/p'=1$.
Note that
\begin{eqnarray*}
 \left | (T_{A_1,A_2} f)(x)\right|
 \le \int_{|y|\le \Omega} |K_{A_1,A_2}(x,y) f(y)|dy \le \int_{|y|\le \Omega} K (x,y) |f(y)|dy.
\end{eqnarray*}
By the dominated convergence theorem, we have
\begin{equation}\label{eq:s2}
  \lim_{A\rightarrow \infty }  \|(T_{A_1,A_2} f - 2\pi\overline{g(x_0)} f )\cdot
   \chi^{}_{\bbR\setminus[-2\Omega,2\Omega]}\|_p=0.
\end{equation}
Now the conclusion follows by combining   (\ref{eq:s1}) and (\ref{eq:s2}).
\end{proof}

We are now ready to give the proof of the main result.

\begin{proof}[Proof of Theorem~\ref{thm:main}]
First, we prove the convergence in $L^p(\bbR)$.

 For any $f\in\LTR$, by (\ref{eq:t2a}), we have
\begin{eqnarray}
(T_{A_1,A_2}f)(x)
&=& \frac{1}{2\pi}
   \int_{\bbR}
     \hat{f}(y)   e^{iyx} dy
     \int_{y-A_2}^{y+A_1}\overline{\hat{g}(\omega)} e^{-ix_0\omega}  d\omega.
  \label{eq:t2}
\end{eqnarray}
Hence
\[
  (T_{A_1,A_2}f) \hat{}\, (y) = h_{A_1,A_2}(y) \hat{f}(y),
\]
where
\[
  h_{A_1,A_2}(y) = \int_{y-A_2}^{y+A_1}\overline{\hat{g}(\omega)} e^{-ix_0\omega}  d\omega.
\]
Obviously, $h_{A_1,A_2}$ is of   bounded variation on $\bbR$ and $V_{h_{A_1,A_2}}\le 2\|\hat g\|_1$.

By Lemma~\ref{Lm:L1} and Proposition~\ref{prop:p4},
$T_{A_1,A_2}$ is a bounded linear  operator on $L^p(\bbR)$ and
\begin{equation}\label{eq:aa1}
    \|T_{A_1,A_2} f\| \le 2C_p  \|\hat g\|_1 \|f\|_p, \qquad \forall f\in L^p(\bbR).
\end{equation}
Fix some $f\in L^p(\bbR)$. For any $\varepsilon>0$, there is some $\tilde f\in C_c^1(\bbR)$ such that
$\|f - \tilde f\|_p < \varepsilon$. By Lemma~\ref{Lm:L2}, we can find some $A_0>0$ such that for any
$A_1, A_2>A_0$,
\[
    \| T_{A_1,A_2} \tilde f -  2\pi\overline{g(x_0)} \tilde f\|_p < \varepsilon.
\]
Consequently,
\begin{eqnarray*}
  \|T_{A_1,A_2} f - 2\pi\overline{g(x_0)}  f\|_p
  &\le& \|T_{A_1,A_2}(f-\tilde f)\|_p +  \| T_{A_1,A_2} \tilde f -  2\pi\overline{g(x_0)} \tilde f\|_p \\
 &&\qquad   +  2\pi|g(x_0)|\cdot \|f - \tilde f\|_p\\
  &\le& (2C_p  \|\hat g\|_1+2\pi|g(x_0)|+ 1)\varepsilon.
\end{eqnarray*}
Hence
\[
    \lim_{A_1,A_2\rightarrow\infty}   \|T_{A_1,A_2} f - 2\pi\overline{g(x_0)}  f\|_p=0, \qquad  \forall f\in L^p(\bbR).
\]

Next we consider the pointwise convergence.
For $A>0$, let $S_A$ be defined by (\ref{eq:sa}).
Then $S_A$ is a bounded linear operator on $L^p(\bbR)$ and
\[
    (S_A f)\hat{} = \hat{f}\cdot \chi^{}_{[-A,A]},\qquad f\in\LTR.
\]
For $f\in L^2\bigcap L^p(\bbR)$, define
\[
  (\tilde T_Af)(x) = \frac{1}{2\pi}
   \int_{\bbR}
 \overline{\hat{g}(\omega)} e^{-ix_0 \omega} (M_{-\omega} S_A M_{\omega}f)(x) d\omega,
\]
where the operator $M_{\omega}$ is defined by
\[
    (M_{\omega}f)(x) = e^{-ix\omega} f(x).
\]
Since $\hat g\in L^1(\bbR)$ and $L^2\bigcap L^p(\bbR)$ is dense in $L^p(\bbR)$,
$\tilde T_A$ can be extended to a bounded linear operator on $L^p(\bbR)$.

On the other hand, for any $f\in L^2\bigcap L^p(\bbR)$, we see from (\ref{eq:t2a}) that
\begin{eqnarray}
  (T_Af)(x)
  &=& \frac{1}{2\pi}
   \int_{\bbR}
 \overline{\hat{g}(\omega)} e^{-ix_0 \omega}  d\omega
     \int_{\omega-A}^{\omega+A}
          \hat{f}(y)   e^{iyx} dy \nonumber \\
  &=& \frac{1}{2\pi}
   \int_{\bbR}
 \overline{\hat{g}(\omega)} e^{-ix_0 \omega}  d\omega
     \int_{-A}^{A}
          \hat{f}(y+\omega)   e^{iyx} e^{ix\omega}dy \nonumber \\
  &=& \frac{1}{2\pi}
   \int_{\bbR}
 \overline{\hat{g}(\omega)} e^{-ix_0 \omega} (M_{-\omega} S_A M_{\omega}f)(x) d\omega \nonumber \\
 &=& (\tilde T_A f)(x). \nonumber
\end{eqnarray}
Using the density of  $L^2\bigcap L^p(\bbR)$ again, we get that $T_A = \tilde T_A$ on $L^p(\bbR)$. Hence
\begin{equation}\label{eq:t21}
     (T_Af)(x) = \frac{1}{2\pi}
   \int_{\bbR}
 \overline{\hat{g}(\omega)} e^{-ix_0 \omega} (M_{-\omega} S_A M_{\omega}f)(x) d\omega,\quad \forall f\in L^p(\bbR).
\end{equation}
It follows that
\[
 \sup_{A> 0} | (T_Af)(x) |
\le \frac{1}{2\pi} \int_{\bbR}|\hat g(\omega)|\cdot   \sup_{A>0}|(M_{-\omega} S_A  M_{\omega}f)(x)|  d\omega.
\]
By   Minkovski's inequality and Proposition~\ref{prop:CH}, we have
\begin{eqnarray}
   \left\|\sup_{A>0} | (T_Af)(x) |\right\|_p
&\le& \frac{1}{2\pi} \int_{\bbR}|\hat g(\omega)|\cdot \left\| \sup_{A> 0}|(M_{-\omega} S_A  M_{\omega}f)(x)|\right\|_p  d\omega
   \nonumber \\
&=& \frac{1}{2\pi} \int_{\bbR}|\hat g(\omega)|\cdot  \left\| \sup_{A>0}|(  S_A  M_{\omega}f)(x)|\right\|_p  d\omega
   \nonumber \\
&\le& \frac{1}{2\pi} \int_{\bbR}|\hat g(\omega)|\cdot  C_p \|M_{\omega}f\|_p  d\omega
   \nonumber \\
&=& \frac{C_p}{2\pi} \|\hat g\|_1 \|f\|_p,\quad \forall f\in L^p(\bbR). \label{eq:ta1}
\end{eqnarray}

Fix some $f\in L^p(\bbR)$.  For any $\varepsilon>0$, we can find some $\tilde f\in C_c^1(\bbR)$ such that
\[
   \| f - \tilde f\|_p <\varepsilon.
\]
By Lemma~\ref{Lm:L2}, we have
\[
 \lim_{A\rightarrow\infty} \| T_A \tilde f - 2\pi\overline{g(x_0)}\tilde f\|_{\infty}=0.
\]
Note that $T_A f$ is continuous on $\bbR$, thanks to Lemma~\ref{Lm:L1}. We have
\[
 \lim_{A\rightarrow\infty}\sup_{x\in\bbR} \left| (T_A \tilde f)(x) - 2\pi\overline{g(x_0)}\tilde f(x)\right|=0.
\]
Hence
\[
   \limsup_{A,A'\rightarrow\infty} |(T_A\tilde f)(x) - (T_{A'}\tilde f)(x)| =0,  \qquad \forall x\in\bbR.
\]
It follows that
\begin{eqnarray*}
&&   \left\|\limsup_{A,A'\rightarrow\infty}\left|(T_A  f)(x) - (T_{A'} f)(x)\right|\right\|_p
    \\
&\le& \left\| \limsup_{A\rightarrow\infty}\left|(T_A (f-\tilde f))(x)\right|\right\|_p
+
   \left\|\limsup_{A,A'\rightarrow\infty} \left|(T_A\tilde f)(x) - (T_{A'}\tilde f)(x)\right|\right\|_p\\
&&\qquad
     +   \left\| \limsup_{A'\rightarrow\infty}\left|(T_{A'} (f-\tilde f))(x)\right|\right\|_p
          \\
&\le&       2 \left\|\sup_{A>0}\left|(T_A(f-\tilde f))(x)\right| \right\|_p
   \\
&\le& \frac{C_p}{\pi} \|\hat g\|_1  \| f - \tilde f\|_p  \qquad \mbox{(using (\ref{eq:ta1}))} \\
&<&\frac{C_p}{\pi} \|\hat g\|_1\cdot \varepsilon.
\end{eqnarray*}
Since $\varepsilon$ is arbitrary, we have
\[
  \left\|\limsup_{A,A'\rightarrow\infty}|(T_A  f)(x) - (T_{A'} f)(x)|\right\|_p=0.
\]
Hence the limit $\lim_{A\rightarrow \infty} (T_Af)(x)$ exists almost everywhere.
Since
$T_A f$ tends to $2\pi\overline{g(x_0)}f$ in $L^p(\bbR)$, we have
\[
  \lim_{A\rightarrow \infty} (T_Af)(x)=2\pi\overline{g(x_0)}f(x),\qquad a.e.
\]
This completes the proof.
\end{proof}

\end{document}